\documentclass[11pt, reqno]{amsart}
\usepackage{amsmath,amssymb,amsfonts,amsthm,graphicx,mathrsfs,mathtools}
\usepackage{array}
\usepackage[margin=1.5in]{geometry}
\usepackage[dvipsnames]{xcolor}
\usepackage{graphicx}
\usepackage[safe]{tipa}
\usepackage{blkarray}
\usepackage{tikz}
\usepackage{array}
\usepackage{tikz}
\usetikzlibrary{matrix,arrows.meta}
\usetikzlibrary{matrix,arrows.meta,calc}
\usepackage{subcaption}
\makeatletter
\@namedef{subjclassname@2020}{%
\textup{2020} Mathematics Subject Classification}
\makeatother
\title[Refined Osserman--Schiffer inequality]{On the refined Osserman--Schiffer inequality in Euclidean space}
\author[J. Lee]{Jaehoon Lee}
\address[]{Jaehoon Lee, Department of Mathematics, The Chinese University of Hong Kong, Shatin, N.T., Hong Kong}
\email{jaehoonlee@cuhk.edu.hk}

\author[E. Yeon]{Eungbeom Yeon}
\address[]{Eungbeom Yeon, Department of Mathematical Sciences, Pusan National University, Busan 46241, Republic of Korea}
\email{ebeom.yeon@pusan.ac.kr}

\allowdisplaybreaks

\begin{document}

\newtheorem{theorem}{theorem}[section]
\newtheorem{thm}[theorem]{Theorem}
\newtheorem{lemma}[theorem]{Lemma}
\newtheorem{cor}[theorem]{Corollary}
\newtheorem{prop}[theorem]{Proposition}
\newtheorem{rmk}[theorem]{Remark}
\newtheorem{Ex}[theorem]{Example}
\newtheorem{conj}[theorem]{Conjecture}
\newtheorem*{mainthm1}{Main Theorem}
\newtheorem*{mainthm2}{Theorem 2}

\numberwithin{equation}{section}
\newcommand{\RNum}[1]{\uppercase\expandafter{\romannumeral #1\relax}}
\newcommand{\R}{\mathbb{R}}
\newcommand{\C}{\mathbb{C}}
\newcommand{\grad}{\nabla}
\newcommand{\laplacian}{\Delta}
\newcommand{\tgamma}{\tilde{\gamma}}
\newcommand{\ttau}{\tilde{\tau}}
\newcommand{\pp}{\Phi}
\newcommand{\tkappa}{\tilde{\kappa}}
\renewcommand{\d}{\textup{d}}
\newlength{\mywidth}
\newcommand\bigfrown[2][\textstyle]{\ensuremath{%
  \array[b]{c}\text{\resizebox{\mywidth}{.7ex}{$#1\frown$}}\\[-1.3ex]#1#2\endarray}}
\newcommand{\arc}[1]{{%
  \setbox9=\hbox{#1}%
  \ooalign{\resizebox{\wd9}{\height}{\texttoptiebar{\phantom{A}}}\cr#1}}}
\newcommand\buildcirclepm[1]{%

  \begin{tikzpicture}[baseline=(X.base), inner sep=-0.1, outer sep=-2]

    \node[draw,circle] (X)  {\footnotesize\raisebox{1ex}{$#1\pm$}};

  \end{tikzpicture}%

}
\newcommand{\introheading}[1]{%
\par\smallskip
\noindent\textbf{#1.}
}

\subjclass[2020]{53A10, 53C42}
\keywords{Osserman--Schiffer inequality, oval conjecture, catenoid, Lagrangian catenoid, Hoffman--Osserman catenoid}

\begin{abstract}
In this paper, we strengthen the Osserman--Schiffer inequality and generalize it to arbitrary codimension, thereby settling a conjecture of Feinberg in a stronger form. Our refinement provides a lower bound for the difference between the two sides of the classical inequality in terms of a nonnegative integral involving the second fundamental form. We show that equality in the classical inequality characterizes concentric planar annuli and rotationally symmetric annuli in catenoids. Furthermore, we characterize the equality cases of the refined inequality under the assumption that equality holds on two distinct sections. Apart from concentric planar annuli, these are precisely circle-foliated annular parts of the surfaces in Hoffman and Osserman's classification of complete doubly connected minimal surfaces in Euclidean space with total curvature $-4\pi$.
\end{abstract}

\maketitle

\section{Introduction}\label{intro}
\setcounter{equation}{0}
The classical Osserman--Schiffer inequality \cite{OssermanSchiffer} provides a sharp convexity estimate for the lengths of the images of concentric circles under conformal harmonic immersions of annuli into $\mathbb{R}^3$. More precisely, the corresponding length function $L(t)$ satisfies $L''(t)\geq L(t)$, where $e^t$ is the radius of the circle in the parameter domain. This estimate has played an important role in the study of isoperimetric inequalities. It has also been used to establish area comparison results for minimal annuli, leading to a variational characterization of the catenoid and nonexistence results for minimal surfaces in a slab \cite{BernsteinBreiner}. In a related direction, a different convexity argument was employed to extend the area comparison to minimal surfaces of more general topology \cite{ChoeDaniel}.

Assuming the oval conjecture of Benguria and Loss \cite{BenguriaLoss}, Bernstein and Breiner proved the following refinement of the Osserman--Schiffer inequality for embedded minimal annuli in a slab.
\begin{prop}[{\cite[Proposition~5.3]{BernsteinBreiner}}]
Let $\Omega$ be the open slab bounded by two parallel planes $P_1=\{x_3=h_1\}$ and $P_2=\{x_3=h_2\}$ in $\mathbb{R}^3$, where $h_1<h_2$. Let $\Sigma\in\mathcal{A}(\Omega)$, and for $h\in (h_1,h_2)$, set $\Sigma_h=\Sigma\cap\{x_3=h\}$. Then
\begin{align*}
\frac{d^2}{dh^2}\mathcal{H}^1(\Sigma_h)\geq \frac{(2\pi)^2}{F_3^2}\mathcal{H}^1(\Sigma_h)+\int_{\Sigma_h}\frac{\left|A(\nu, E_2)\right|^2}{\left|\nabla_{\Sigma} x_3\right|^2}\,d\mathcal{H}^1,
\end{align*}
where $A$ is the second fundamental form of $\Sigma$, and $\{\nu, E_2\}$ is a global orthonormal frame on $\Sigma$ with $\nu$ parallel to $\nabla_{\Sigma}x_3$.
\end{prop}
Here, $\mathcal{A}(\Omega)$ denotes the class of smooth embedded minimal annuli $\Sigma\subset\Omega$ with $\overline{\Sigma}\setminus\Sigma\subset P_1\cup P_2$. The quantity $F_3$ denotes the vertical flux of $\Sigma$.

The above inequality can be reformulated in terms of conformal harmonic immersions of annuli. In this paper, we show that this refined inequality holds for every conformal harmonic immersion of an annulus into $\mathbb{R}^n$, for all $n\geq 3$, without any embeddedness or slab assumptions and independently of the oval conjecture.

To state our main result, let $X: A_{r_1,r_2}\to\mathbb{R}^n$, $n\geq 3$, be a smooth conformal harmonic immersion, where $A_{r_1,r_2}:=\{\zeta\in\mathbb{C}\mid r_1<|\zeta|<r_2\}$. Using logarithmic polar coordinates $\zeta=e^{t+i\theta}$, conformality gives
\begin{align*}
\langle X_t, X_\theta\rangle=0,\quad |X_t|=|X_\theta|=:\rho>0,
\end{align*}
where $\rho=\rho(t,\theta)$ is the conformal factor. Define the unit tangent and conormal vectors by
\begin{align*}
\tau:=\frac{X_\theta}{\rho},\quad \nu:=\frac{X_t}{\rho}.
\end{align*}

We establish the following inequality.

\begin{thm}\label{RefinedOS}
Let 
\begin{align*}
L(t):=\int_0^{2\pi}\rho(t,\theta)\,d\theta,\quad t\in (\log r_1, \log r_2).
\end{align*}
Then
\begin{align*}
L''(t)\geq L(t)+\int_0^{2\pi}\rho^3\left|A(\nu, \tau)\right|^2\,d\theta,
\end{align*}
where $A$ denotes the second fundamental form of the immersion $X$.
\end{thm}

Our theorem generalizes the Osserman--Schiffer inequality to higher codimension and settles a conjecture of Feinberg \cite[Section~6, Question~1]{Feinberg} in a stronger form. In higher codimension, Feinberg established weaker convexity estimates of the form $L''(t)>\lambda L(t)$ for some $0<\lambda<1$ and conjectured that the sharp Osserman--Schiffer inequality would still hold. This conjecture was subsequently recorded as an open problem by Osserman in his survey on the isoperimetric inequality \cite[Section~4, Question~4]{Osserman}. 

The following theorem shows that the characterization of equality cases in \cite[Theorem~1]{OssermanSchiffer} also extends naturally from $\R^3$ to $\R^n$.

\begin{thm}\label{thm stability}
The equality
\begin{align*}
L''(t_0)=L(t_0)
\end{align*}
holds for some $t_0\in(\log r_1,\log r_2)$ if and only if $X$ is a one-to-one parametrization of either a concentric planar annulus or a rotationally symmetric annulus in a catenoid. In both cases, $L''=L$ holds throughout $(\log r_1,\log r_2)$.
\end{thm}

In the next theorem, we further analyze the equality cases of the refined inequality $L''\geq L+\mathcal R$, where
\begin{align}\label{eq:mixed-integral}
\mathcal R(t)=\int_0^{2\pi}\rho^3|A(\nu,\tau)|^2\,d\theta.
\end{align}
To state the theorem, we first introduce the following family of minimal surfaces:
\begin{equation}\label{eq:catenoid-models}
\begin{gathered}
X_{a,b,c}(t,\theta)=\bigl(a\cosh t\cos\theta,\ a\cosh t\sin\theta,\ b\sinh t\cos\theta,\ b\sinh t\sin\theta,\ ct\bigr),\\
a^2=b^2+c^2,\quad a>0,\quad b,c\ge0,\quad e^{t+i\theta}\in \C\setminus\{0\}.
\end{gathered}
\end{equation} 

Write $X_{a,b,c}(t,\theta)=(x_1,\ldots,x_5)$. If $b=0$, then $a=c$, and \eqref{eq:catenoid-models} gives the standard parametrization of a catenoid in $\R^3$. If $c=0$, then $a=b$, and the orthogonal change of coordinates
\begin{align}\label{eq:lagrangian-coordinates}
z=\frac{x_1+x_3+i(x_2+x_4)}{\sqrt2},\quad
w=\frac{x_1-x_3-i(x_2-x_4)}{\sqrt2}
\end{align}
gives
\begin{align*}
zw=\frac{a^2}{2}.
\end{align*}
Thus, in this case, the image becomes the Lagrangian catenoid in $\R^4$. 

Suppose now that $b,c>0$ and introduce the complex coordinate $\zeta=e^{t+i\theta}$. After an orthogonal transformation, the parametrization \eqref{eq:catenoid-models} can be written as
\begin{align}\label{eq:ho-representation}
\operatorname{Re}\left(-\frac{c^2}{4\kappa}\zeta-\frac{\kappa}{\zeta},\
i\frac{c^2}{4\kappa}\zeta-\frac{i\kappa}{\zeta},\
c\log\zeta,\ f\zeta,\ if\zeta\right),
\end{align}
where
$$
\kappa=\frac{\sqrt{a^2+b^2}}2,\quad f=\frac{ab}{\sqrt{a^2+b^2}}.
$$
Up to rigid motions and scaling, the surface in \eqref{eq:ho-representation} is uniquely determined by the ratio $\frac{\kappa f}{c^2}$ (see \cite[Remark~3.3]{LeeYeon}). 

The above observations show that \eqref{eq:catenoid-models} represents the entire circle-foliated family in the Hoffman--Osserman classification \cite[Proposition~6.6]{HO}, up to rigid motions and scaling. Among these examples, we refer to those that fully lie in $\mathbb{R}^5$ as Hoffman--Osserman catenoids.
 
The following theorem shows that the family parametrized by (\ref{eq:catenoid-models}) gives all possible nonplanar equality cases of the refined inequality.
\begin{thm}\label{thm2 stability}
The equality
\begin{align*}
L''(t_i)=L(t_i)+\mathcal R(t_i),\quad i=1,2
\end{align*}
holds for two distinct $t_1,t_2\in(\log r_1,\log r_2)$ if and only if $X$ is a one-to-one parametrization of either a concentric planar annulus or a circle-foliated annular part of a surface in the Hoffman--Osserman classification. The circle-foliated family in this classification consists of catenoids, Lagrangian catenoids, and Hoffman--Osserman catenoids. In all these cases, $L''=L+\mathcal R$ holds throughout $(\log r_1,\log r_2)$. 
\end{thm}


\section*{Acknowledgements}
JL was supported substantially by research grants from the Research Grants Council of the Hong Kong Special Administrative Region, China [Project No.: CUHK 14306624]. EY was supported by National Research Foundation of Korea NRF-2022R1C1C2013384. 

\introheading{On the use of AI}\enspace The authors formulated the problem and initiated this project without using AI. We first tried to solve the problem using the Weierstrass representation, following the approach of Osserman and Schiffer. When we encountered difficulties with this approach, ChatGPT suggested applying the Wirtinger inequality to the function $q$. This approach directly leads to the refined inequality in $\mathbb{R}^3$. We used this observation to extend the argument to $\mathbb{R}^n$ for all $n\geq 3$, and this extension process was done without AI assistance. Without using AI, we also found that the only possible nonplanar examples are annular parts of the circle-foliated surfaces in the Hoffman--Osserman classification. ChatGPT was also used to improve the English text.


\section{Proof of Theorem \ref{RefinedOS}}\label{sec2}
\setcounter{equation}{0}
We begin by recalling the classical total curvature inequality, originally established by Fenchel \cite{Fenchel} in $\mathbb{R}^3$ and later extended to $\mathbb{R}^n$ by Borsuk \cite{Borsuk}. For each $t\in (\log r_1, \log r_2)$, it gives
\begin{align}\label{totalcurv}
K(t):=\int_0^{2\pi}\omega(t,\theta)\,d\theta\geq 2\pi,
\end{align}
where $\omega:=|\tau_\theta|$. Equality holds if and only if the corresponding curve is a convex plane curve. See also \cite[Theorem~3.4]{Milnor}.

Define the functions $\psi$ and $q$ on $(\log r_1, \log r_2)\times [0,2\pi]$ by
\begin{align}\label{QQ}
\psi(t,\theta):= \frac{2\pi}{K(t)}\int_0^\theta \omega(t,s)\,ds,\quad
q(t,\theta):=\sqrt{\rho(t,\theta)}\,e^{\frac{i\psi(t,\theta)}{2}}.
\end{align}
Since $\psi(t,0)=0$ and $\psi(t,2\pi)=2\pi$, we have
\begin{align*}
q(t,2\pi)=-q(t,0).
\end{align*}
Thus, for each fixed $t$, the function $q(t,\cdot)$ admits an antiperiodic extension to $\mathbb{R}$ satisfying
\begin{align*}
q(t,\theta+2\pi)=-q(t,\theta).
\end{align*} 

Applying the Wirtinger inequality for antiperiodic functions, we obtain
\begin{align}\label{Wir}
\int_0^{2\pi}|q_\theta(t,\theta)|^2\,d\theta\geq\frac{1}{4}\int_0^{2\pi}|q(t,\theta)|^2\,d\theta,
\end{align}
where equality holds if and only if $q(t,\theta)=c_0(t)e^{\frac{i\theta}{2}}+c_{-1}(t)e^{-\frac{i\theta}{2}}$ for some $c_0(t),c_{-1}(t)\in\mathbb{C}$.

\begin{lemma}\label{lem21}
For all $t\in (\log r_1, \log r_2)$, we have
\begin{align*}
\int_0^{2\pi}\frac{|X_{\theta\theta}(t,\theta)|^2}{\rho(t,\theta)}\,d\theta\geq\int_0^{2\pi}\rho(t,\theta)\,d\theta.
\end{align*}
\end{lemma}
\begin{proof}
Differentiating (\ref{QQ}) with respect to $\theta$ gives
\begin{align*}
q_\theta=\left(\frac{\rho_\theta}{2\sqrt{\rho}}+i\frac{\pi}{K}\sqrt{\rho}\omega\right)e^{\frac{i\psi}{2}},\quad 4|q_\theta|^2=\frac{\rho_\theta^2}{\rho}+\frac{4\pi^2}{K^2}\rho\omega^2.
\end{align*}
Hence, by (\ref{Wir}), 
\begin{align}\label{ineq1}
\int_0^{2\pi}\frac{\rho_\theta^2}{\rho}+\frac{4\pi^2}{K^2}\rho\omega^2\,d\theta=\int_0^{2\pi}4|q_\theta|^2\,d\theta\geq\int_0^{2\pi}|q|^2\,d\theta=\int_0^{2\pi}\rho\,d\theta.
\end{align}

On the other hand, differentiating $X_\theta=\rho\tau$ yields
\begin{align*}
X_{\theta\theta}=\rho_\theta\tau+\rho\tau_\theta.
\end{align*}
Since $|\tau|=1$ and $\omega=|\tau_\theta|$, we have
\begin{align*}
\frac{|X_{\theta\theta}|^2}{\rho}=\frac{\rho_\theta^2}{\rho}+\rho\omega^2.
\end{align*}
Moreover, $K\geq2\pi$ by (\ref{totalcurv}), and hence
\begin{align}\label{ineq2}
\int_0^{2\pi}\frac{|X_{\theta\theta}|^2}{\rho}\,d\theta=\int_0^{2\pi}\frac{\rho_\theta^2}{\rho}+\rho\omega^2\,d\theta\geq\int_0^{2\pi}\frac{\rho_\theta^2}{\rho}+\frac{4\pi^2}{K^2}\rho\omega^2\,d\theta.
\end{align}

Combining (\ref{ineq1}) and (\ref{ineq2}), we obtain
\begin{align*}
\int_0^{2\pi}\frac{|X_{\theta\theta}|^2}{\rho}\,d\theta\geq\int_0^{2\pi}\rho\,d\theta,
\end{align*}
which completes the proof.
\end{proof}

\begin{lemma}\label{comp}
We have
\begin{align*}
\rho_{tt}=\frac{\rho_\theta^2}{\rho}+\rho^3|A(\nu,\tau)|^2-\frac{\langle X_\theta, X_{\theta\theta\theta}\rangle}{\rho}.
\end{align*}
\end{lemma}
\begin{proof}
Differentiating $|X_t|^2=\rho^2$ with respect to $\theta$ and $|X_\theta|^2=\rho^2$ with respect to $t$, we obtain
\begin{align*}
\langle X_{t\theta}, \nu\rangle=\rho_\theta,\quad \langle X_{t\theta}, \tau\rangle=\rho_t.
\end{align*}
Let $(X_{t\theta})^\perp$ denote the orthogonal projection of $X_{t\theta}$ onto the normal space of the immersion. Since $(X_{t\theta})^\perp=\rho^2A(\nu,\tau)$, we have
\begin{align*}
|X_{t\theta}|^2&=\langle X_{t\theta}, \nu\rangle^2+\langle X_{t\theta}, \tau\rangle^2+|(X_{t\theta})^\perp|^2\\
&=\rho_\theta^2+\rho_t^2+\rho^4|A(\nu,\tau)|^2.
\end{align*}

Differentiating $|X_\theta|^2=\rho^2$ twice with respect to $t$ gives
\begin{align*}
\rho\rho_{tt}+\rho_t^2=|X_{t\theta}|^2+\langle X_\theta, X_{\theta tt}\rangle.
\end{align*}
Substituting the expression for $|X_{t\theta}|^2$ yields
\begin{align*}
\rho_{tt}=\frac{\rho_\theta^2}{\rho}+\rho^3|A(\nu,\tau)|^2+\frac{\langle X_\theta, X_{\theta tt}\rangle}{\rho}.
\end{align*}

Since $X$ is harmonic, $X_{tt}+X_{\theta\theta}=0$, and hence
\begin{align*}
\langle X_\theta, X_{\theta tt}\rangle=-\langle X_\theta, X_{\theta\theta\theta}\rangle.
\end{align*}
This completes the proof.
\end{proof}

We are now ready to prove the theorem.
\begin{proof}[Proof of Theorem \ref{RefinedOS}]
By Lemma \ref{comp}, we have
\begin{align*}
L''(t)=\int_0^{2\pi}\rho_{tt}\,d\theta=\int_0^{2\pi}\frac{\rho_\theta^2}{\rho}+\rho^3|A(\nu,\tau)|^2-\frac{\langle X_\theta, X_{\theta\theta\theta}\rangle}{\rho}\,d\theta.
\end{align*}
Integration by parts, using periodicity, gives
\begin{align*}
-\int_0^{2\pi}\frac{\langle X_\theta, X_{\theta\theta\theta}\rangle}{\rho}\,d\theta&=\int_0^{2\pi}\frac{|X_{\theta\theta}|^2}{\rho}\,d\theta-\int_0^{2\pi}\frac{\rho_\theta\langle X_\theta, X_{\theta\theta}\rangle}{\rho^2}\,d\theta\\
&=\int_0^{2\pi}\frac{|X_{\theta\theta}|^2}{\rho}\,d\theta-\int_0^{2\pi}\frac{\rho_\theta^2}{\rho}\,d\theta.
\end{align*}
Hence,
\begin{align*}
L''(t)=\int_0^{2\pi}\frac{|X_{\theta\theta}|^2}{\rho}\,d\theta+\int_0^{2\pi}\rho^3|A(\nu,\tau)|^2\,d\theta.
\end{align*}

By Lemma \ref{lem21}, 
\begin{align*}
\int_0^{2\pi}\frac{|X_{\theta\theta}|^2}{\rho}\,d\theta\geq\int_0^{2\pi}\rho\,d\theta=L(t).
\end{align*}
Consequently,
\begin{align*}
L''(t)\geq L(t)+\int_0^{2\pi}\rho^3|A(\nu,\tau)|^2\,d\theta,
\end{align*}
which completes the proof.
\end{proof}
\section{Equality cases}\label{sec3}

Throughout this section, let $X:A_{r_1,r_2}\to\R^n$ be as in Theorem~\ref{RefinedOS}. Thus,
\begin{align*}
|X_t|=|X_\theta|=\rho>0,\quad \langle X_t,X_\theta\rangle=0,\quad X_{tt}+X_{\theta\theta}=0.
\end{align*}
The proof of Theorem~\ref{RefinedOS} gives
\begin{align}\label{eq:length-identity}
L''(t)=\int_0^{2\pi}\frac{|X_{\theta\theta}|^2}{\rho}\,d\theta+\mathcal R(t).
\end{align}

We first determine the equality case in Lemma~\ref{lem21}, whose proof applies to every regular $C^2$ closed curve. If $\gamma$ parametrizes a circle in the form
\begin{align}\label{eq:circle}
\gamma(\theta)=p+R(e_1\cos\theta+e_2\sin\theta),
\end{align}
where $p\in\R^n$, $R>0$, and $e_1,e_2$ are orthonormal vectors in $\R^n$, then equality holds in Lemma~\ref{lem21}. The following lemma shows that the converse also holds.

\begin{lemma}\label{lem:curve-equality}
Let $\gamma:\mathbb{S}^1\to\R^n$ be a regular $C^2$ closed curve, with angular parameter $\theta$. The equality
$$
\int_0^{2\pi}\frac{|\gamma_{\theta\theta}|^2}{|\gamma_\theta|}\,d\theta=\int_0^{2\pi}|\gamma_\theta|\,d\theta
$$
holds if and only if $\gamma$ has the form \eqref{eq:circle}.
\end{lemma}

\begin{proof}
Suppose that
$$
\int_0^{2\pi}\frac{|\gamma_{\theta\theta}|^2}{|\gamma_\theta|}\,d\theta=\int_0^{2\pi}|\gamma_\theta|\,d\theta.
$$
Put
\begin{align*}
\rho_\gamma=|\gamma_\theta|,\quad T=\frac{\gamma_\theta}{\rho_\gamma},\quad \omega=|T_\theta|,\quad K=\int_0^{2\pi}\omega\,d\theta,
\end{align*}
and define $\psi$ and $q$ as in \eqref{QQ}, with $\rho$ replaced by $\rho_\gamma$. The argument of Lemma~\ref{lem21}, applied to $\gamma$, shows that equality holds in both \eqref{ineq1} and \eqref{ineq2}. Since $\rho_\gamma>0$ and $K\ge2\pi$ by Fenchel's theorem, we have $\int_0^{2\pi}\rho_\gamma\omega^2\,d\theta>0$. Equality in \eqref{ineq2} therefore implies that $K=2\pi$. 

By the equality case of Fenchel's theorem, $\gamma$ is a convex plane curve. Identify the plane containing $\gamma$ with $\C$, choosing the orientation so that the signed curvature of $\gamma$ is nonnegative. Since $K=2\pi$, \eqref{QQ} gives $\psi_\theta=\omega$ and $\psi(0)=0$. The chosen orientation gives $T_\theta=i\omega T$, so $T(\theta)=T(0)e^{i\psi(\theta)}$. Therefore,
\begin{align}\label{eq:tangent-square}
\gamma_\theta=T(0)q^2.
\end{align}

Equality in \eqref{Wir} gives
\begin{align}\label{eq:q-equality}
q(\theta)=c_0e^{i\theta/2}+c_{-1}e^{-i\theta/2}
\end{align}
for some $c_0,c_{-1}\in\C$. Since $\gamma$ is closed, \eqref{eq:tangent-square} and \eqref{eq:q-equality} imply
\begin{align*}
0=\int_0^{2\pi}\gamma_\theta\,d\theta=4\pi T(0)c_0c_{-1}.
\end{align*}

At least one of $c_0,c_{-1}$ is nonzero because $|q|^2=\rho_\gamma>0$. Moreover, \eqref{QQ} gives
\begin{align}\label{eq:q-orientation}
\operatorname{Im}(q_\theta\overline q)=\frac12\rho_\gamma\omega\ge0.
\end{align}
If $c_0=0$, then the left-hand side of \eqref{eq:q-orientation} becomes $-\frac{|c_{-1}|^2}{2}<0$. Hence $c_{-1}=0$ and $c_0\ne0$. By \eqref{eq:tangent-square}, $\gamma_\theta=T(0)c_0^2e^{i\theta}$. Integrating this identity with respect to $\theta$ gives \eqref{eq:circle} with $R=|c_0|^2$.
\end{proof}

\begin{proof}[Proof of Theorem~\ref{thm stability}]
Since $(X_{t\theta})^\perp=\rho^2A(\nu,\tau)$, \eqref{eq:mixed-integral} gives
$$
\mathcal R(t)=\int_0^{2\pi}\frac{|(X_{t\theta})^\perp|^2}{\rho}\,d\theta.
$$
Suppose $L''(t_0)=L(t_0)$. Since $\mathcal R\ge0$, we have
\begin{align*}
\mathcal R(t_0)=0,\quad L(t_0)=\int_0^{2\pi}\frac{|X_{\theta\theta}(t_0,\theta)|^2}{\rho(t_0,\theta)}\,d\theta.
\end{align*}
By Lemma~\ref{lem:curve-equality}, a suitable translation gives
\begin{align*}
X(t_0,\theta)=Re_r(\theta),\quad e_r=e_1\cos\theta+e_2\sin\theta,\quad e_\theta=\partial_\theta e_r,
\end{align*}
for some $R>0$ and orthonormal vectors $e_1,e_2\in\R^n$. 

At $t=t_0$, we have $X_t=R\nu$. Since $\rho>0$ and the integrand in \eqref{eq:mixed-integral} is continuous and nonnegative, $\mathcal R(t_0)=0$ implies that $A(\nu,\tau)=0$ along the section $\theta\mapsto X(t_0,\theta)$. The identity $X_{t\theta}=R\nu_\theta$ then shows that $\nu_\theta$ is tangent to the surface. Together with $\langle\nu_\theta,\nu\rangle=0$, this implies that $\nu_\theta$ is a scalar multiple of $e_\theta$. Along the section $t=t_0$, set $\beta(\theta)=\langle\nu_\theta,e_\theta\rangle$. Then $\nu_\theta=\beta e_\theta$. Differentiating $\langle\nu,e_\theta\rangle=0$ with respect to $\theta$ and using $\partial_\theta e_\theta=-e_r$, we obtain
\begin{align*}
\beta=\langle\nu,e_r\rangle.
\end{align*}
It follows that $\beta_\theta=\langle\nu_\theta,e_r\rangle+\langle\nu,e_\theta\rangle=0$, so $\beta$ is constant. Integrating $\nu_\theta=\beta e_\theta$ gives
\begin{align*}
\nu=\beta e_r+\eta
\end{align*}
for a constant vector $\eta\in\R^n$. Since $\langle\nu,e_\theta\rangle=0$ for all $\theta$ and $|\nu|=1$, we have
\begin{align}\label{eq:ordinary-constraints}
\eta\perp e_1,e_2,\quad \beta^2+|\eta|^2=1.
\end{align}

Now, we consider the following harmonic map
\begin{align}\label{eq:ordinary-model}
\widetilde X(t,\theta)=R(\cosh s+\beta\sinh s)e_r(\theta)+Rs \eta,\quad s=t-t_0.
\end{align}
At $t=t_0$, $\widetilde X$ and $X$ have the same value and $t$-derivative. To show that $X=\widetilde X$, put $Y=X-\widetilde X$ and define
\begin{align*}
\widehat Y_j(t)=\frac1{2\pi}\int_0^{2\pi}Y(t,\theta)e^{-ij\theta}\,d\theta,\quad j\in\mathbb{Z}.
\end{align*}
Since $Y$ is harmonic and periodic in $\theta$, integration by parts gives
\begin{align*}
\widehat Y_j''=j^2\widehat Y_j,\quad \widehat Y_j(t_0)=\widehat Y_j'(t_0)=0.
\end{align*}
Uniqueness for the initial value problem gives $\widehat Y_j\equiv0$. For each $t$, all Fourier coefficients of $Y(t,\cdot)$ vanish, so $Y(t,\cdot)=0$ in $L^2(\mathbb S^1;\R^n)$. By continuity, $Y\equiv0$, and hence $X=\widetilde X$.

If $\eta=0$, then \eqref{eq:ordinary-constraints} gives $\beta=\pm1$. Thus $X=Re^{\pm s}e_r(\theta)$, which parametrizes a planar annulus.
If $\eta\ne0$, then $|\beta|<1$ by \eqref{eq:ordinary-constraints}. Choose $s_0\in\R$ such that $\tanh s_0=\beta$, and put
\begin{align*}
a=R|\eta|,\quad e_3=\frac{\eta}{|\eta|}.
\end{align*}
Then $\cosh s+\beta\sinh s=|\eta|\cosh(s+s_0)$, and \eqref{eq:ordinary-model} gives
\begin{align*}
X(t,\theta)+Rs_0\eta=a\bigl(\cosh(s+s_0)e_r(\theta)+(s+s_0)e_3\bigr).
\end{align*}
Thus $X$ parametrizes a rotationally symmetric annulus contained in a catenoid. 

We remark that $X$ is one-to-one in both cases. One can easily check that $L''=L$ and $\mathcal R\equiv0$ for the planar annuli and annular parts of catenoids described in the theorem, which completes the proof.
\end{proof}

\begin{lemma}\label{lem:two-circles}
Suppose that for distinct $t_1,t_2$ the curves
$$
\theta\mapsto X(t_i,\theta),\quad i=1,2,
$$
are one-to-one constant-speed parametrizations of circles. Then $X$ extends to a conformal minimal immersion on $\C\setminus\{0\}$. For each $t\in\R$, the section $\theta\mapsto X(t,\theta)$ is a one-to-one constant-speed parametrization of a circle.
\end{lemma}

\begin{proof}
For $j\in\mathbb{Z}$, define
\begin{align*}
\widehat X_j(t)=\frac1{2\pi}\int_0^{2\pi}X(t,\theta)e^{-ij\theta}\,d\theta.
\end{align*}
By harmonicity,
\begin{equation}\label{eq:fourier-ode}
\widehat X_j''=j^2\widehat X_j.
\end{equation}
By \eqref{eq:circle}, each component of $X(t,\theta)$ at $t=t_1,t_2$ is a linear combination of $1$, $\cos\theta$, and $\sin\theta$. Hence
$$
\widehat X_j(t)=\frac{1}{2\pi}\int_0^{2\pi}X(t,\theta)e^{-ij\theta}\,d\theta=0,\quad t\in\{t_1,t_2\},\quad |j|\ge2.
$$
For such $j$, substituting $\widehat X_j(t_1)=0$ into the general solution of \eqref{eq:fourier-ode} gives
$$
\widehat X_j(t)=\frac{\widehat X_j'(t_1)}{|j|}\sinh\bigl(|j|(t-t_1)\bigr).
$$
Since $t_1\ne t_2$, substituting $t=t_2$ and using $\widehat X_j(t_2)=0$ gives $\widehat X_j'(t_1)=0$, and hence $\widehat X_j\equiv0$. Applying \eqref{eq:fourier-ode} to the remaining Fourier coefficients $\widehat X_0$ and $\widehat X_{\pm1}$, we obtain
\begin{align*}
X(t,\theta)=p_0+(t-t_1)d+U(t)\cos\theta+V(t)\sin\theta,\quad U''=U,\quad V''=V
\end{align*}
for some constant vectors $p_0,d\in\R^n$ and $\R^n$-valued functions $U,V$ of $t$. 

Writing $X(t_1,\cdot)$ in the form \eqref{eq:circle} gives $U(t_1)=Re_1$ and $V(t_1)=Re_2$. Comparing Fourier coefficients in $\langle X_t,X_\theta\rangle=0$ at $t=t_1$ gives
\begin{align}\label{eq:initial-orthogonality}
d\perp e_1,e_2,\quad U'(t_1)\perp e_2,\quad V'(t_1)\perp e_1,\quad \langle U'(t_1),e_1\rangle=\langle V'(t_1),e_2\rangle=\alpha
\end{align}
for some $\alpha\in\R$. Similarly, $|X_t|^2=|X_\theta|^2=R^2$ at $t=t_1$ gives
\begin{align}\label{eq:initial-lengths}
d\perp U'(t_1),V'(t_1),\quad |U'(t_1)|=|V'(t_1)|,\quad U'(t_1)\perp V'(t_1),\quad |d|^2+|U'(t_1)|^2=R^2.
\end{align}

Let $u=U'(t_1)-\alpha e_1$ and $v=V'(t_1)-\alpha e_2$ denote the orthogonal projections of $U'(t_1)$ and $V'(t_1)$, respectively, onto $\operatorname{span}\{e_1,e_2\}^{\perp}$. Equations \eqref{eq:initial-orthogonality} and \eqref{eq:initial-lengths} give
\begin{align}\label{eq:two-circle-constraints}
&u,v\perp e_1,e_2,\quad |u|=|v|,\quad \langle u,v\rangle=0,\nonumber\\
&d\perp e_1,e_2,u,v,\quad R^2=\alpha^2+|u|^2+|d|^2.
\end{align}
Solving $U''=U$ and $V''=V$ with these initial data gives
\begin{align}\label{eq:two-circle-form}
X(t,\theta)&=p_0+(t-t_1)d+\sinh(t-t_1)(u\cos\theta+v\sin\theta)\nonumber\\
&\quad+\bigl(R\cosh(t-t_1)+\alpha\sinh(t-t_1)\bigr)(e_1\cos\theta+e_2\sin\theta).
\end{align}

Formula \eqref{eq:two-circle-form} defines a harmonic map for all $t\in\R$. By \eqref{eq:two-circle-constraints}, $U(t)$ and $V(t)$ are orthogonal, with common length $r(t)$. Thus the radius of the section $X(t,\cdot)$ satisfies
\begin{align}\label{eq:radius}
r(t)^2=\bigl(R\cosh(t-t_1)+\alpha\sinh(t-t_1)\bigr)^2+|u|^2\sinh^2(t-t_1).
\end{align}
Differentiating \eqref{eq:two-circle-form} and using \eqref{eq:two-circle-constraints}, we obtain $\langle X_t,X_\theta\rangle=0$ and $|X_t|^2=|X_\theta|^2=r^2$. 

Since $|\alpha|\le R$ and $R>0$, we have $R\cosh(t-t_1)+\alpha\sinh(t-t_1)>0$ for all $t\in\R$. Thus every map of the form \eqref{eq:two-circle-form} with coefficients satisfying \eqref{eq:two-circle-constraints} is a conformal minimal immersion on $\C\setminus\{0\}$, and each section has the form \eqref{eq:circle}.
\end{proof}

\begin{proof}[Proof of Theorem~\ref{thm2 stability}]
By \eqref{eq:length-identity} and Lemma~\ref{lem:curve-equality}, equality in Theorem~\ref{RefinedOS} at a fixed $t_0\in(\log r_1,\log r_2)$ holds if and only if
\begin{align}\label{eq:one-circle}
X(t_0,\theta)=p+R(e_1\cos\theta+e_2\sin\theta)
\end{align}
for some $p\in\R^n$, $R>0$, and orthonormal vectors $e_1,e_2\in\R^n$. Suppose equality in Theorem~\ref{RefinedOS} holds on two distinct sections $t=t_1$ and $t=t_2$. By \eqref{eq:two-circle-form}, we can write
$$
X(t,\theta)=p_0+(t-t_1)d+U(t)\cos\theta+V(t)\sin\theta,
$$
for some $\R^n$-valued functions $U,V$ satisfying
$$
U(t)\perp V(t),\quad |U(t)|=|V(t)|=r(t)>0
$$
for every $t\in(\log r_1,\log r_2)$. Thus each section is a circle parametrized as in \eqref{eq:one-circle}, so equality in Theorem~\ref{RefinedOS} holds throughout this interval.

With the same notation as in the proof of Lemma~\ref{lem:two-circles}, set $c=|d|$. If $u=0$ and $c=0$, then \eqref{eq:two-circle-constraints} gives $v=d=0$ and $\alpha=\pm R$. After a suitable translation, we obtain
\begin{align*}
X(t,\theta)=Re^{\pm(t-t_1)}(e_1\cos\theta+e_2\sin\theta),\quad R>0.
\end{align*}
This is a one-to-one parametrization of a concentric planar annulus. 

Otherwise, $|u|^2+c^2>0$, so $|\alpha|<R$ by \eqref{eq:two-circle-constraints}. By \eqref{eq:radius}, we obtain
\begin{align*}
4r(t)^2=\bigl((R+\alpha)^2+|u|^2\bigr)e^{2(t-t_1)}+\bigl((R-\alpha)^2+|u|^2\bigr)e^{-2(t-t_1)}+2c^2.
\end{align*}
Both coefficients of the exponential functions are positive, so $r$ has a unique minimum on $\R$. By Lemma~\ref{lem:two-circles}, $X$ extends to all $t\in\R$, and every section of the extension is a one-to-one constant-speed parametrization of a circle. After translating $t$ so that $r$ attains its minimum at $t=0$, we may take $t_1=0$ in \eqref{eq:two-circle-form}. Equation \eqref{eq:radius} then gives $\alpha=r'(0)=0$.

By \eqref{eq:two-circle-constraints}, the vectors $U(0)$, $V(0)$, $U'(0)$, $V'(0)$, and $d$ are mutually orthogonal and have lengths $a$, $a$, $b$, $b$, and $c$, respectively, where $a^2=b^2+c^2$ and $a>0$. Solving $U''=U$ and $V''=V$ with these initial data gives \eqref{eq:catenoid-models} up to a rigid motion. Thus the image lies in an affine subspace of $\R^n$ of dimension at most five. Its affine span has dimension three if $b=0$, four if $c=0$, and five if $b,c>0$. To see this, evaluate \eqref{eq:catenoid-models} at $\theta=0$ and $\theta=\pi/2$ and use the linear independence of $1$, $t$, $\cosh t$, and $\sinh t$.

To see that $X$ is one-to-one, first suppose $c>0$. The fifth coordinate determines $t$, and the first two coordinates determine $\theta$. If $c=0$, then \eqref{eq:lagrangian-coordinates} gives $z=\frac{a}{\sqrt{2}} e^{t+i\theta}$, whose modulus determines $t$ and whose argument determines $\theta$ modulo $2\pi$. Thus each restriction is a one-to-one parametrization of an annulus formed by circles.

Conversely, for concentric planar annuli and annular restrictions of \eqref{eq:catenoid-models}, every one-to-one conformal parametrization has sections of the form \eqref{eq:one-circle}, so equality holds for all $t\in(\log r_1,\log r_2)$.
\end{proof}

\begin{rmk}\normalfont
Let $\gamma$ be a one-to-one constant-speed parametrization of a horizontal circle $C$ in $\R^3$. Given any real analytic unit normal field along $C$, the Bj\"orling construction yields a minimal annulus near $C$ for which equality holds on this section \cite[Section~2]{MeeksWeber}. If the vertical component of the prescribed normal field is nonconstant along $C$, then Theorem~\ref{thm2 stability} implies that equality holds on no other section.
\end{rmk}

\end{document}